\documentclass[12pt]{amsart}
\usepackage{amsmath,amsfonts,amssymb,amscd}
\usepackage{latexsym}

\title [Diagonal coinvariants and DAHA]
{Diagonal coinvariants and Double affine Hecke algebras}
\author[Ivan Cherednik]{Ivan Cherednik $^\dag$} 
\date{May 15, 2003}

\thanks{$^\dag$ \ Partially supported by NSF grant 
DMS-0200276}

\address[I. Cherednik]{Department of Mathematics, UNC 
Chapel Hill, North Carolina 27599, USA\\
chered@email.unc.edu}

\newcommand{\Z}{{\mathbb Z}}
\newcommand{\Q}{{\mathbb Q}}
\newcommand{\N}{{\mathbb N}}
\newcommand{\C}{{\mathbb C}}
\newcommand{\R}{{\mathbb R}}

\def\HH{\mbox{${\mathcal H}$\kern-5.2pt${\mathcal H}$}}

\newtheorem{theorem}{Theorem}[section]

\newtheorem{definition}[theorem]{Definition}
\newtheorem{lemma}[theorem]{Lemma}

\newtheorem{theorem }{Theorem}[section]
\newtheorem{proposition }[theorem]{Proposition}
\newtheorem{definition }[theorem]{Definition}
\newtheorem{lemma }[theorem]{Lemma}
\newtheorem{corollary }[theorem]{Corollary}
\newtheorem{notation }[theorem]{Notation}
\newtheorem{remark }[theorem]{Remark}
\newtheorem{example }[theorem]{Example}

\newtheorem{ theorem}{Theorem}[section]
\newtheorem{ proposition}[theorem]{Proposition}
\newtheorem{ definition}[theorem]{Definition}
\newtheorem{ lemma}[theorem]{Lemma}
\newtheorem{ corollary}[theorem]{Corollary}
\newtheorem{ notation}[theorem]{Notation}
\newtheorem{ remark}[theorem]{Remark}
\newtheorem{ example}[theorem]{Example}

 \newcommand{\rem}{{\bf Comment.\ }}

\def\for{\  \hbox{ for } \ }
\def\iif{ \ \hbox{ if } \ }

\def\where{\  \hbox{ where } \ }
\def\and{\  \hbox{ and } \ }

\def\equal{\stackrel{\,\mathbf{def}}{= \kern-3pt =}}

\def\om{\omega}

\def\al{\alpha}

\def\de{\delta}
\def\De{\Delta}

\def\si{\sigma}

\def\Ga{\Gamma}
\def\ze{\zeta}

\def\kapp{\hbox{\bf \ae}}

\def\vep{\varepsilon}

\def\vth{{\vartheta}}

\def\tal{\tilde{\alpha}}

\def\tV{\tilde{V}}

\def\tGa{\tilde{\Gamma}}

\def\tw{\widetilde w}
\def\tW{\widetilde W}

\def\tV{\tilde V}
\def\tz{\tilde z}
\def\tb{\tilde b}

\def\tR{\tilde R}

\def\hw{\widehat{w}}
\def\hW{\widehat{W}}

\def\hv{\hat{v}}

\def\F{\mathbf{F}}

\def\0{\mathbf{0}}
\def\H{\mathbf{H}}

\def\çF{\mathcal{F}}

\def\d{\mathcal{D}}
\def\p{\mathcal{P}}

\def\h{\mathcal{H}}

\def\g{\mathcal{G}}

\def\i{\mathcal{I}}

\def\lan{\langle}

\def\ran{\rangle}

\def\lng{\hbox{\tiny lng}}
\def\sht{\hbox{\tiny sht}}

\def\inv{\hbox{\tiny inv}}

\newcommand{\sgn}{\mbox{sgn}}

\def\HH{\mathfrak{H}}

\def\BB{\mathfrak{B}}

\def\HH{\hbox{${\mathcal H}$\kern-5.2pt${\mathcal H}$}}

\font\smm=msbm10 at 12pt 
\def\symbol#1{\hbox{\smm #1}}
\def\lsmash{{\symbol n}}

\def\#{\sharp}

\begin{document}
\maketitle
{\small
\tableofcontents
}

\vfil

It was conjectured by Haiman \cite{Hai} that the space of
diagonal coinvariants for a root system $R$ of rank $n$
has a "natural" quotient of dimension $(1+h)^n$ for
the Coxeter number $h$. This space 
is the quotient $\C[x,y]/(\C[x,y]\C[x,y]^W_o)$ 
for the algebra of polynomials $\C[x,y]$ with the diagonal action
of the Weyl group on $x\in \C^n \ni y$ and the ideal 
$\C[x,y]^W_o \subset \C[x,y]^W$ of the $W$-invariant 
polynomials without the constant term.
In \cite{Go}, such a quotient was constructed.
It appeared to be the graded object of the perfect
module (in the terminology of
\cite{C12}) of the rational double affine Hecke algebra
for the simplest nontrivial $k=-1-1/h.$

Generally, the perfect modules are defined as irreducible
self-dual spherical representations
of DAHA with a projective action of the $PSL_2(\Z).$
In the $q$-case, the semisimplicity is added.
At roots of unity, they generalize the Verlinde algebras.
Gordon gives an explicit description of the above 
module as a quotient of the space of double polynomials
considered as a representation
of the rational DAHA induced from the sign-character of 
the nonnaffine Weyl group $W.$
 
In \cite{C12}, perfect modules appear
naturally as quotients of the algebra of single
polynomials. Using the double polynomials has some
advantages. For instance, the self-duality becomes obvious.
We note that the perfect modules are always quotients of 
the space of double polynomials, but the corresponding 
kernels are expected to be reasonably simple 
only for $k\in -1/h-\Z_+.$   
\medskip

We extend Gordon's description to the $q$-case,
establishing its {\em direct} connection with a fundamental
fact that the Weyl algebra of rank $n$ 
(a noncommutative $n$-torus) has
a unique irreducible representation provided that 
the center element $q$ is a primitive $N$-th
root of unity and the generators are cyclic of order $N.$
Its dimension is $N^n,$ which matches the Haiman number as  
$N=1+h.$ We {\em deduce} our theorem from this fact
and, as a corollary, obtain a new, entirely algebraic, proof
of Gordon's theorem using the Lusztig-type 
isomorphism acting from
the general DAHA to its rational degeneration.

Gordon's demonstration was based on the results 
due to Opdam-- Rouquier
(see \cite{GGOR}) on the monodromy of the KZ-connection 
from \cite{Ch0},\cite{Ch1}
in relation to the representation theory of the rational DAHA.
The technique of Lusztig's isomorphisms is actually 
of the same origin. These isomorphisms are closely connected with the monodromy
of the {\em affine} KZ-connection.

\medskip
The structure of the paper is as follows.
We start 
with the general definitions and the construction of
the Lusztig-type isomorphisms (for reduced root systems). 
Then we switch to the case of generic $q$ and
$k=-1-1/h.$  Following \cite{C7},\cite{C12}, we make $q$ a
root of unity. In \cite{C7}, this method was used to 
deform the
Verlinde algebras. Finally, we obtain our theorem and 
reprove (the main part of) the theorem from \cite{Go}. 

\smallskip
Lusztig's isomorphisms are important by themselves. 
Under minor restrictions, they establish
an equivalence of the categories of finite dimensional
representations of DAHA when $q$ is not a root of unity
and those of its rational degeneration.
They can be applied to infinite dimensional 
representations as well, however,
generally speaking, the theory gets analytic.  
\smallskip

The last section of the
paper contains the definitions of the following two
new objects, the {\em universal} double affine Hecke algebra 
and the corresponding universal Dunkl operators 
acting in the {\em noncommutative} polynomials in terms of 
two sets of variables $X$ and $Y.$ 
Upon the reduction to the commutative polynomials,
these operators are directly connected with the main theorem 
and have other applications. We note that the 
universal DAHA satisfies a "noncommutative" variant 
of the PBW--theorem.

\medskip
The definition of the
universal DAHA is $X\leftrightarrow Y$--symmetric as well
as its homomorphism to the DAHA. 
The $X\leftrightarrow Y$--duality of
the DAHA was deduced in \cite{C15} from the topological 
interpretation of the double affine ("elliptic") 
braid group. 
A direct proof of the DAHA-duality is not difficult too
(see \cite{Ma3}). 
The universal DAHA, to be more
exact the corresponding braid group, can be used
to simplify the proof.
It has something in common with the method from \cite{IS}.

Concerning the elliptic braid group and 
its topological interpretation,
there is a connection with the construction due to v.d.Lek,
although the orbifold fundumental group was used in \cite{C15}
insead of removing the ramification divisor 
in his construction.
This connection was mentioned in \cite{C15} 
and is discussed in 
more detail in \cite{Io}. In the case of $GL_n,$ the braid group from \cite{C15} is essentially due to Birman and Scott.

\medskip
The author is thankful to P.~Etingof for 
important remarks and discussions. 
The paper was completed at Institut 
de Math\'ematiques de Luminy.
I am grateful for the invitation. 
I also thank the referee for a thorough report.  

\medskip
\vskip 0.2cm
\section{Double affine Hecke algebras}
\setcounter{equation}{0}

Let $R=\{\al\}\subset \R^n$ be a root system of type $A,B,...,F,G$
with respect to a euclidean form $(z,z')$ on $\R^n 
\ni z,z'$,
$W$ the Weyl group  generated by the reflections $s_\al$,
$R_{+}$ the set of positive  roots, 
corresponding to (fixed) simple roots  
roots $\al_1,...,\al_n,$ 
$\Ga$ the Dynkin diagram  
with $\{\al_i, 1 \le i \le n\}$ as the vertices, 
$R^\vee=\{\al^\vee =2\al/(\al,\al)\}$ the dual root system,
\begin{align}
& Q=\oplus^n_{i=1}\Z \al_i \subset P=\oplus^n_{i=1}\Z \om_i, 
\notag \end{align}
where $\{\om_i\}$ are fundamental weights:
$ (\om_i,\al_j^\vee)=\de_{ij}$ for the 
simple coroots $\al_i^\vee.$

The form will be normalized
by the condition  $(\al,\al)=2$ for the 
{\em short} roots. This normalization coincides with that
from the tables in \cite{Bo} for $A,C,D,E,G.$
Hence $\nu_\al\equal (\al,\al)/2$ can be  $1,2$ or $3.$  
Sometimes we write $\nu_{\lng}$ for long roots ($\nu_{\sht}= 1$).
Let  $\vth\in R^\vee $ be the maximal positive 
{\em coroot} (it is maximal short in $R$),
$\rho= (1/2)\sum_{\al\in R_+} \al \ =
\sum_i \om_i.$


\vskip 0.2cm
{\bf Affine roots.}
The vectors $\ \tal=[\al,\nu_\al j] \in 
\R^n\times \R \subset \R^{n+1}$ 
for $\al \in R, j \in \Z $ form the 
{\em affine root system} 
$\tR \supset R$ ($z\in \R^n$ are identified with $ [z,0]$).  
We add $\al_0 \equal [-\vth,1]$ to the simple
roots.
The set $\tR$ of positive roots is
$R_+\cup \{[\al,\nu_\al j],\ \al\in R, \ j > 0\}$.
Let $\tal^\vee=\tal/\nu_\al,$ so $\al_0^\vee=\al_0.$

The Dynkin diagram $\Ga$ of $R$  
is completed by $\al_0$ (by $-\vth$ to be more
exact). The notation is $\tGa$. It is the
completed Dynkin diagram for $R^\vee$ from \cite{Bo}
with the arrows reversed.

The set of
the indices of the images of $\al_0$ by all 
the automorphisms of $\tGa$ will be denoted by $O$ 
($O=\{0\} \for E_8,F_4,G_2$). Let $O'={r\in O, r\neq 0}$.
The elements $\om_r$ for $r\in O'$ are the so-called minuscule
weights: $(\om_r,\al^\vee)\le 1$ for
$\al \in R_+$.

Given $\tal=[\al,\nu_\al j]\in \tR,  \ b \in P$, let  
\begin{align}
&s_{\tal}(\tz)\ =\  \tz-(z,\al^\vee)\tal,\ 
\ b(\tz)\ =\ [z,\ze-(z,b)]
\label{ondonx}
\end{align}
for $\tz=[z,\ze] \in \R^{n+1}$.

The {\it affine Weyl group} $\tW$ is generated by all $s_{\tal}$. 
One can take
the simple reflections $s_i=s_{\al_i}\ (0 \le i \le n)$ 
as its generators and introduce the corresponding notion of the  
length. This group is
the semidirect product $W\lsmash Q$ of 
its subgroups $W$ and the lattice $Q$, 
where $\al\in Q$ is identified with
$ s_{\al}s_{[\al,\nu_{\al}]}=$ 
$s_{[-\al,\nu_\al]}s_{\al}$ for
$\al\in R.$

The {\it extended Weyl group} $ \hW$ generated by $W\and P$
is isomorphic to $W\lsmash P$:
\begin{align}
&(wb)([z,\ze])\ =\ [w(z),\ze-(z,b)] \for w\in W,\, b\in P.
\label{ondthrx}
\end{align}

Given $b\in P_+$, let $w^b_0$ be the longest element
in the subgroup $W_0^{b}\subset W$ of the elements
preserving $b$. This subgroup is generated by simple 
reflections. We set
\begin{align}
&u_{b} = w_0w^b_0  \in  W,\ \pi_{b} =
b( u_{b})^{-1}
\ \in \ \hW, \  u_i= u_{\om_i},\pi_i=\pi_{\om_i},
\label{wox}
\end{align}
where $w_0$ is the longest element in $W,$
$1\le i\le n.$

The elements $\pi_r\equal\pi_{\om_r}, r \in O'$ and
$\pi_0=\hbox{id}$ leave $\tGa$ invariant 
and form a group denoted by $\Pi$, 
 which is isomorphic to $P/Q$ by the natural 
projection $\{\om_r \mapsto \pi_r\}$. As to $\{ u_r\}$,
they preserve the set $\{-\vth,\al_i, i>0\}$.
The relations $\pi_r(\al_0)= \al_r= ( u_r)^{-1}(-\vth)$ 
distinguish the
indices $r \in O'$. Moreover,
\begin{align}
& \hW  = \Pi \lsmash \tW, \where
  \pi_rs_i\pi_r^{-1}  =  s_j \iif \pi_r(\al_i)=\al_j,\ 
 0\le j\le n.
\end{align}

Setting
$\hw = \pi_r\tw \in \hW,\ \pi_r\in \Pi,\, \tw\in \tW,$
the {\it length} $l(\hw)$ 
is by definition the length of the reduced decomposition 
$\tw = s_{i_l}...s_{i_2} s_{i_1} $
in terms of the simple reflections 
$s_i, 0\le i\le n.$ 

\medskip
{\bf DAHA.}
By  $m,$ we denote the least natural number 
such that  $(P,P)=(1/m)\Z.$  Thus
$m=2 \for D_{2k},\ m=1 \for B_{2k}, C_{k},$
otherwise $m=|\Pi|$.

The double affine Hecke algebra depends 
on the parameters 
$q, t_\nu,\, \nu\in \{\nu_\al\}.$ The definition ring is 
$\Q_{q,t}\equal$
$\Q[q^{\pm 1/m},t^{\pm 1/2}]$ formed by the
polynomials in terms of $q^{\pm 1/m}$ and  
$\{t_\nu^{\pm 1/2} \}.$
We set
\begin{align}
&   t_{\tal} = t_{\al}=t_{\nu_\al},\ t_i = t_{\al_i},\ 
q_{\tal}=q^{\nu_\al},\ q_i=q^{\nu_{\al_i}},\notag\\ 
&\where \tal=[\al,\nu_\al j] \in \tR,\ 0\le i\le n.
\label{taljx}
\end{align}

It will be convenient to use the parameters
$\{k_\nu\}$ together with  $\{t_\nu \},$ setting
$$
t_\al=t_\nu=q_\al^{k_\nu} \for \nu=\nu_\al, \and
\rho_k=(1/2)\sum_{\al>0} k_\al \al.
$$

For pairwise commutative $X_1,\ldots,X_n,$    
\begin{align}
& X_{\tb}\ =\ \prod_{i=1}^nX_i^{l_i} q^{ j} 
\iif \tb=[b,j],\ \hw(X_{\tb})\ =\ X_{\hw(\tb)}.
\label{Xdex}
\\  
&\hbox{where\ } b=\sum_{i=1}^n l_i \om_i\in P,\ j \in 
\frac{1}{ m}\Z,\ \hw\in \hW.
\notag \end{align}
Later $Y_{\tb}=Y_b q^{-j}$ will be needed. Note the
opposite sign of $j$. We set $(\tb,c)=(b,c).$

We will also use that $\pi_r^{-1}$ is $\pi_{r^*}$ and
$u_r^{-1}$ is $u_{r^*}$ 
for $r^*\in O\ ,$  $u_r=\pi_r^{-1}\om_r.$
The reflection $^*$ is 
induced by an involution of the nonaffine Dynkin diagram
$\Gamma.$

\begin{definition}
The  double  affine Hecke algebra $\HH\ $
is generated over $ \Q_{ q,t}$ by 
the elements $\{ T_i,\ 0\le i\le n\}$, 
pairwise commutative $\{X_b, \ b\in P\}$ satisfying 
(\ref{Xdex}),
and the group $\Pi,$ where the following relations are imposed:

(o)\ \  $ (T_i-t_i^{1/2})(T_i+t_i^{-1/2})\ =\ 
0,\ 0\ \le\ i\ \le\ n$;

(i)\ \ \ $ T_iT_jT_i...\ =\ T_jT_iT_j...,\ m_{ij}$ 
factors on each side;

(ii)\ \   $ \pi_rT_i\pi_r^{-1}\ =\ T_j \iif 
\pi_r(\al_i)=\al_j$; 

(iii)\  $T_iX_b T_i\ =\ X_b X_{\al_i}^{-1} \iif 
(b,\al^\vee_i)=1,\
0 \le i\le  n$;

(iv)\ $T_iX_b\ =\ X_b T_i$ if $(b,\al^\vee_i)=0 
\for 0 \le i\le  n$;

(v)\ \ $\pi_rX_b \pi_r^{-1}\ =\ X_{\pi_r(b)}\ =\ 
X_{ u^{-1}_r(b)}
 q^{(\om_{r^*},b)},\  r\in O'$.
\label{doublex}
\end{definition}
\qed

Given $\tw \in \tW, r\in O,\ $ the product
\begin{align}
&T_{\pi_r\tw}\equal \pi_r\prod_{k=1}^l T_{i_k},\where 
\tw=\prod_{k=1}^l s_{i_k},
l=l(\tw),
\label{Twx}
\end{align}
does not depend on the choice of the reduced decomposition
(because $\{T\}$ satisfy the same ``braid'' relations 
as $\{s\}$ do).
Moreover,
\begin{align}
&T_{\hv}T_{\hw}\ =\ T_{\hv\hw}\  \hbox{ whenever}\ 
 l(\hv\hw)=l(\hv)+l(\hw) \for
\hv,\hw \in \hW. \label{TTx}
\end{align}
In particular, we arrive at the pairwise 
commutative elements 
\begin{align}
& Y_{b}\ =\  \prod_{i=1}^nY_i^{l_i} \iif  
b=\sum_{i=1}^n l_i\om_i\in P,\where  
 Y_i\equal T_{\om_i},
\label{Ybx}
\end{align}
satisfying the relations
\begin{align}
&T^{-1}_iY_b T^{-1}_i\ =\ Y_b Y_{\al_i}^{-1} \iif 
(b,\al^\vee_i)=1,
\notag\\ 
& T_iY_b\ =\ Y_b T_i \iif (b,\al^\vee_i)=0,
 \ 1 \le i\le  n.
\end{align}
For arbitrary nonzero $q,t,$ any element $H \in \HH\ $  
has a unique decomposition in the form
\begin{align}
&H =\sum_{w\in W }\,  g_{w}\, f_w\, T_w,\ 
g_{w} \in \Q_{q,t}[X],\ f_{w} \in \Q_{q,t}[Y],  
\label{hatdecx}
\end{align}
and five more analogous decompositions corresponding 
to the other orderings
of $\{T,X,Y\}.$ It makes the polynomial representation
(to be defined next)
the $\HH\,$-module induced from the one dimensional
representation $T_i\mapsto t_i^{1/2},\,$ $Y_i\mapsto Y_i^{1/2}$
of the affine Hecke subalgebra $\h_Y=\lan T,Y\ran.$

These and below statements are from \cite{C2}.

One may also use the 
{\em intermediate subalgebras} of \HH\
with $P$ replaced by any lattice $B\ni b$ between $Q$ and $P$
for $X_b$ and $Y_b$ (see \cite{C12}). Respectively, $\Pi$ is 
changed to the preimage of $B/Q$ in $\Pi.$ Generally, 
there can be two different
lattices $B_x$ and $B_y$ for $X$ and $Y.$
The $m\in \N$ from the definition of $\Q_{q,t}$
has to be the least such that $m(B_x,B_y)\subset \Z.$

Note that $\HH\, ,$ its degenerations, and the 
corresponding polynomial representations
are actually defined over $\Z$ extended by the parameters of DAHA.
We will use its $\Z$-structure a couple of times in the paper
(the modular reduction), but prefer to stick to $\Q.$
The Lusztig isomorphisms require $\Q.$

\vskip 0.2cm
{\bf Demazure-Lusztig operators.} They
are defined as follows: 
\begin{align}
&T_i\  = \  t_i ^{1/2} s_i\ +\ 
(t_i^{1/2}-t_i^{-1/2})(X_{\al_i}-1)^{-1}(s_i-1),
\ 0\le i\le n,
\label{Demazx}
\end{align}
and obviously preserve $\Q[q,t^{\pm 1/2}][X]$.
We note that only the formula for $T_0$ involves $q$: 
\begin{align}
&T_0\  =  t_0^{1/2}s_0\ +\ (t_0^{1/2}-t_0^{-1/2})
( q X_{\vth}^{-1} -1)^{-1}(s_0-1),\notag\\ 
&\where
s_0(X_b)\ =\ X_bX_{\vth}^{-(b,\vth)}
 q^{(b,\vth)},\ 
\al_0=[-\vth,1].
\end{align}

The map sending $ T_j$ to the formula from
(\ref{Demazx}), $\ X_b \mapsto X_b$ 
(see (\ref{Xdex})),
$\pi_r\mapsto \pi_r$ induces a 
$ \Q_{ q,t}$-linear 
homomorphism from $\HH\, $ to the algebra of linear endomorphisms 
of $\Q_{ q,t}[X]$.
This $\HH\,$-module, which will be called the
{\em polynomial representation}, 
is faithful 
and remains faithful when   $q,t$ take  
any nonzero complex values assuming that
$q$ is not a root of unity. 

The images of the $Y_b$ are called the 
{\em difference Dunkl operators}. To be more exact,
they must be called trigonometric-difference Dunkl operators,
because there are also 
rational-difference Dunkl operators.

\vskip 0.2cm
{\bf Automorphisms.}
Assuming that $B_x=B_y,$ 
the following maps can be uniquely extended to 
automorphisms of
$\HH\ $(see \cite{C4},\cite{C12}):
\begin{align}
\vep:\ X_i& \mapsto Y_i,\   Y_i \mapsto X_i,\   
 T_i \mapsto T_i^{-1}\,(i\ge 1),\,
t_\nu \mapsto t_\nu^{-1},\,
 q\mapsto  q^{-1},
\label{vepx}
\\ 
 \tau_+:  X_b \mapsto X_b, \ &Y_r \mapsto 
X_rY_r q^{-\frac{(\om_r,\om_r)}{2}},\ 
T_i\mapsto T_i\,(i\ge 1),\ \ t_\nu \mapsto t_\nu,\
 q\mapsto  q,
\notag\\ 
\tau_+:\ &Y_\vth \mapsto q^{-1}\,X_\vth T_0^{-1} 
T_{s_\vth},\, T_0\mapsto q^{-1}\,X_\vth T_0^{-1}, 
\and
\label{taux}\\
\tau_-\ & \equal  \vep\tau_+\vep,\and  
\si\equal \tau_+\tau_-^{-1}\tau_+ =
\tau_-^{-1}\tau_+\tau_-^{-1}= \vep\si^{-1}\vep, 
\label{tauminax}
\end{align}
where $r\in O'.$
In the definition of $\tau_\pm$ and $\si,$ 
we need to add $q^{\pm 1/(2m)}$ to
$\Q_{q,t}.$ 
Here the quadratic relation (o) from Definition 
\ref{doublex} may be omitted. Only the group relations matter. 
The elements $\tau_\pm$ generate the projective $PSL(2,\Z),$
which is isomorphic to the braid group $B_3$ due to Steinberg.

\vskip 0.2cm
{\bf Intertwining operators.}
The {\em $Y$-intertwiners} (see \cite{C1})
are introduced as follows:
\begin{align}
&\Phi_i\ =\  
T_i + (t_i^{1/2}-t_i^{-1/2})
(Y_{\al_i}^{-1}-1)^{-1} \for 1\le i\le n,
\notag\\ 
&\Phi_0\ =\ 
X_\vth T_{s_\vth} - (t_0^{1/2}-t_0^{-1/2})(Y_0 -1)^{-1},
\ Y_0=Y_{\al_0}\equal q^{-1}Y_{\vth}^{-1},
\notag\\
& G_i=\Phi_i (\phi_i)^{-1},\   
\phi_i= t_i^{1/2} + 
(t_i^{1/2} -t_i^{-1/2})(Y_{\al_i}^{-1}-1)^{-1}.
\label{Phix}
\end{align}

Actually these formulas are the $\vep$-images of the formulas
for the $X$-intertwiners, which are a straightforward generalization
of those in the affine Hecke theory.

They belong to $\HH\ $ extended by 
the rational functions in terms of $\{Y\}$. The $G$ are called
the {\em normalized
intertwiners}. The elements
$$
G_i,\ P_r\equal\ X_r T_{u_r^{-1}},\  0\le i\le n,\ r\in O',
$$ 
satisfy the same relations
as $\{s_i,\pi_r\}$ do, so the map 
\begin{align} 
\hw\mapsto G_{\hw}\ =\ P_r G_{i_l}\cdots G_{i_1},
\where \hw=\pi_r s_{i_l}\cdots s_{i_1}\in \hW,
\label{Phiprodx}
\end{align} 
is a  well defined homomorphism  from $\hW.$

The intertwining property is 
$$
G_{\hw} Y_b G_{\hw}^{-1}=Y_{\hw(b)} 
\where Y_{[b,j]}\equal Y_b q^{-j}.
$$
The $P_1$ in the case of $GL$ is due to Knop and Sahi.

As to $\Phi_i$, they
satisfy the  homogeneous Coxeter relations 
and those with $\Pi_r.$ So we may set
$\Phi_{\hw} =$ $P_r \Phi_{i_l}\cdots \Phi_{i_1}$ for the reduced
decompositions. They intertwine $Y$ as well.

The formulas for $\Phi_i$ when $1\le i\le n$ are well known
in the theory of affine Hecke algebras. 
The affine intertwiners are  
the raising operators for the Macdonald nonsymmetric polynomials,
serve the Harish-Chandra -- Opdam spherical transform, and are the 
key tool in the theory of semisimple representations of DAHA.

\vskip 0.2cm
\section {Degenerate DAHA} \label{SECTDEGEN}
\setcounter{equation}{0}

Recall that $m(P,P)\in \Z$ or $m(B,B)\in \Z$ if $B$ is used,
$$ k_i=k_{\al_i},\ k_0=k_{\sht},\ \nu_\al=(\al,\al)/2\in \{1,2,3\}.
$$
We set $\Q_k\equal\Q[k_\al].$ If the integral coefficients
are needed, we take $\Z_k\equal\Z[k_\al,1/m]$ as the definition ring.

The {\em degenerate (graded) double affine Hecke algebra} $\HH' $ is 
the span of 
the group algebra $\Q_k [\hW]$ and the pairwise commutative
\begin{align}
&y_{\tb}\equal
\sum ^n_{i=1}(b,\al_i^\vee)y_i +u \for 
\tb=[b,u]\in P\times\Z,
\notag \end{align}
satisfying  the following relations:
\begin{align}
&s_j y_{b}-y_{s_j(b)}s_j\ =\ - k_j(b,\al_j ),\ 
(b,\al_0)\equal -(b,\vth),
\notag\\ 
&\pi_r y_{\tb}\ =\ y_{\pi_r(\tb)}\pi_r \for 
  0\le j\le n, \  r\in O.
\label{sukax}
\end{align}

\rem
Without $s_0$ and $\pi_r$, we arrive at the defining relations
of the graded affine Hecke algebra from  \cite{L}.
The algebra $\HH'$ has two natural
polynomial representations via the differential-trigonometric and
difference-rational Dunkl operators. There is also the third one,
a representation in terms of infinite differential-trigonometric
Dunkl operators, which leads to 
differential-elliptic $W$-invariant operators generalizing those
due to Ol\-sha\-net\-sky- Pe\-re\-lo\-mov. 
See, e.g., \cite{C1}.
We will need here only  
the (most known) differential-trigonometric polynomial
representations.
\qed

Let us establish the connection with the general DAHA. We set
$$
q=\exp(\mathfrak{v}),\ t_j=q_i^{k_i}=q^{\nu_{\al_i}k_i}, 
\ Y_b=\exp(-\mathfrak{v} y_b),\ \mathfrak{v}\in \C.
$$
Using $\vep$ from (\ref{vepx}),
the algebra $\HH\,$ is generated by $Y_b,\,$ $T_i$ for $1\le i\le n,$ and 
$$
\vep(T_0)=X_{\vth}T_{s_{\vth}},\ 
\vep(\pi_r)=X_rT_{u_r^{-1}},\ r\in O'.
$$
It is straightforward to see that the relations
(\ref{sukax}) for $y_b,\,s_i (i>0),\,$ 
$s_0,\,\pi_r$ are the leading coefficients
of the $\mathfrak{v}$-expansions of the general relations
for this system of generators.
Thus $\HH'$ is $\HH\ $ in the limit $\mathfrak{v}\to 0.$

When calculating the limits of the
$Y_b$ in the polynomial representation,
the "trigonometric" derivatives of $\Q[X]$ appear: 
$$
\partial_{a}(X_{b})\ =\ (a ,b)X_{b},\ a,b\in P,\ \, 
w(\partial _{b})=
\partial_{w(b)}, \  w\in W.
$$ 

The $Y_b$ result in the {\em trigonometric Dunkl operators} 
\begin{align}
\d_{b} \equal\ 
&\partial_b + \sum_{\al\in R_+} \frac{ k_{\al}(b,\al)}{
(1-X_{\al}^{-1}) }
\bigl( 1-s_{\al} \bigr)- (\rho_{k},b).
\label{dunkx}
\end{align}
They act on the Laurent polynomials $f\in \Q_k[X],$ 
are pairwise commutative, and $y_{[b,u]}=\d_b+u$
satisfy (\ref{sukax}) for 
the following action of the group $\hW$:
$$ w^x(f)=w(f)\for w\in W,\ b^x(f)=X_b f \for b\in P.$$
For instance, $s_0^x(f)=X_\vth s_\vth(f),\ \pi_r^x(f)=X_r u_r^{-1}(f).$

Degenerating $\{ \Phi\}$, one gets the intertwiners of 
$\HH'\ $:
\begin{align}
&\Phi'_i= 
s_i + \frac{\nu_i k_i }{ y_{\al_i} },\ 0\le i\le n,\ 
\bigl(\,\Phi'_0= 
X_\vth s_\vth + \frac{k_0 }{ 1-y_{\vth}} \hbox{\ in\ } 
\Q_k[X]\,\bigr),\notag\\ 
&P_r' = \pi_r,\ \bigl(\,P_r'=X_r u_r^{-1} \hbox{\ in\ } 
\Q_k[X]\, \bigr),\ r\in O'.
\label{Phiprimex}
\end{align}
The operator $P_1'$ in the case of $GL$ (it is of 
infinite order) plays the key role in \cite{KnS}.

Recall that the general normalized intertwiners are
$$
G_i=\Phi_i\phi_i^{-1}, \ 
\phi_i=t^{1/2}+(t_i^{1/2}-t^{1/2})(Y_{\al_i}^{-1}-1)^{-1}.
$$
Their limits are
$$G'_i=\Phi'_i(\phi'_i)^{-1},\ 
\phi'_i=1+\frac{\nu_i k_i }{ y_{\al_i} }.$$
They satisfy the unitarity condition
$(G_i')^2=1,$ and the products $G'_{\hw}$ can be defined for any
decompositions of $\hw.$ 
One has: 
$$G'_{\hw}\, y_b\, (G'_{\hw})^{-1}\ =\ y_{\hw(b)}.$$

Equating 
$$G_i=G'_i \for 0\le i\le n,\ P_r=P'_r \for r\in O, 
$$ 
we come to the formulas for $T_i\, (0\le i\le n),\, $ 
$X_r\, (r\in O')$ in terms of 
$s_i, y_b,$ $Y_b=\exp(-\mathfrak{v}y_b).$ 

These formulas determine 
the {\em Lusztig homomorphism} $\kapp'$ from $\HH\, $ to  
the completion $\Z_{k,q,t}\HH'[[\mathfrak{v}y_b]]$ for 
$\Z_{k,q,t}\equal\Z_k\Z_{q,t}.$ 
See, e.g., \cite{C1}.

For instance,
$X_r\in \HH$ becomes $\pi_r T_{u_r^{-1}}^{-1}$ in $\HH',$ 
where the $T$-factor has to be further expressed in terms of $s,y.$
In the degenerate polynomial representation, $\kapp'(X_r)$ acts 
as $X_r (\kapp'(T_{u_r^{-1}})u_r)^{-1},$ not as the 
straightforward multiplication
by $X_r.$  They coincide only 
in the limit $\mathfrak{v}\to 0,$ when $T_w$ become $w.$

Upon the $\mathfrak{v}$-completion,
we get an isomorphism 
$$
\kapp':\Q_k[[\mathfrak{v}]]\otimes\HH\to 
\Q_k[[\mathfrak{v}]]\otimes\HH'.
$$

We will use the
notation $(d,[\al,j])=j.$ For instance,
$(b+d,\al_0)=1-(b,\vth).$ 

Treating $\mathfrak{v}$ 
as a nonzero number,
an arbitrary
$\HH'$-module $V'$ which is a union
of finite dimensional $Y$-modules has a natural structure
of an $\HH$-module provided that we have 
\begin{align}
& q^{(\al_i,\xi+d)}=t_i \Rightarrow 
(\al_i,\xi+d) = \nu_i k_i,\label{lusisox}
\\
& q^{(\al_i,\xi+d)}=\ 1 \Rightarrow 
(\al_i,\xi+d) = 0, \where\notag \\
& 0\le i\le n,\  
y_b(v')=(b,\xi)v' \for \xi\in \C^n,\ 0\neq v'\in V'.\notag  
\end{align}

For the modules of this type,
the map $\kapp'$ is over the ring $\Q_{k,q,t}$ 
extended by $(\al,\xi+d),\ $ $q^{(\al,\xi+d)}$ for $\al\in R$ 
and $y$-eigenvalues $\xi.$  Moreover, we need to localize
by $(1- q^{(\al,\xi+d)})\neq 0$ 
and by $(\al,\xi+d)\neq 0.$
Upon such extension and localization, $\kapp'$  is defined over
$\Z_{k,q,t}$ 
if the module is $y$-semisimple. If there are nontrivial
Jordan blocks, then the formulas will contain factorials in the
denominators. 

For instance, let $\i'[\xi]$ be the $\HH'$-module  induced from the
one-dimensional $y$-module $y_b(v)=(b,\xi)v.$ Assuming that
$q$ is not a root of unity, the mapping $\kapp'$
supplies it with a structure of $\HH$-module if 
$$
q^{(\al,\xi)+\nu_\al j}=t_\al \hbox{\ implies\ }
(\al,\xi)+\nu_\al j=\nu_\al k_\al
$$
for every $\al\in R, j\in \Z,$
and the corresponding implications hold for $t$ replaced by $1.$ 
This means that
\begin{align} 
&(\al,\xi)-\nu_\al k_\al,\, (\al,\xi)\, \not\in \, 
\nu_\al\Z\, +\, \frac{2\pi\imath}{\mathfrak{v}}
(\Z\setminus \{0\})\hbox{\ for\ all\ }\al\in R.
\label{intercond}
\end{align}
Generalizing, $\kapp'$ is well defined for any 
$\HH'$-module generated by its
$y$-eigenvectors with the $y$-eigenvalues $\xi$ satisfying this
condition, assuming that $\mathfrak{v}\not\in \pi\imath\Q.$

\rem
Actually there are at least four different variants of $\kapp'$ because the
normalization factors $\phi,\phi'$ may be associated with different 
one dimensional characters
of the affine Hecke algebra $\lan T,Y\ran$
and its degeneration. There is also a
possibility to multiply the normalized intertwiners by the characters of $\hW$
before equating. Note that if we  
divide the intertwiners $\Phi$ and/or $\Phi'$
by $\phi,\phi'$
on the left in the definition of $G,G',$ it corresponds to
switching from $T_i\mapsto t_i$ to
the character $T_i \mapsto -t_i^{-1/2}$ together with
the multiplication by the sign-character of $\tW.$
In the paper, we will use only $\kapp'$ introduced above. 
\qed

\vskip 0.2cm
{\bf Rational degeneration.} The limit to the Dunkl operators is as follows.
We set $X_b=e^{\mathfrak{w} x_b},\ d_b(x_c)=(b,c),$ so the above derivatives
$\partial_b$ become $\partial_b=(1/\mathfrak{w})d_b.$  
In the limit $\mathfrak{w}\to 0,$
$\mathfrak{w}\d_b$ tends to  
\begin{align}
D_{b} \equal
&d_b +
\sum_{\al\in R_+} \frac{ k_{\al}(b,\al)}{
x_\al }
\bigl( 1-s_{\al} \bigr).
\label{dunkorigx}
\end{align}
These operators are pairwise commutative and satisfy
the cross-relations
\begin{align}
&D_bx_c-x_cD_b=(b,c)+\sum_{\al>0} k_\al (b,\al)(c,\al^\vee)s_\al,
\for b,c\in P.
\label{duncrossx}
\end{align}
These relations, 
the commutativity of $D,$ the commutativity of $x,$ and the $W$-equivariance
$$ w\, x_b\, w^{-1}\,=\,x_{w(b)},\ w\, D_b\, w^{-1}\,=\,D_b\ 
\for b\in P,\, w\in W,
$$
are the defining relations of the {\em rational DAHA} $\HH''.$ 

The references are \cite{CM} (the case of $A_1$) and
\cite{EG}, however the key part of the definition is the
commutativity of $D_b$ due to Dunkl \cite{Du}. 
The Dunkl operators and the operators of multiplication by the $x_b$ 
form the
{\em polynomial representation} of $\HH'',$ which is faithful. It
readily justifies the PBW-theorem for $\HH''.$

Note that in contrast to the $q,t$-setting, the definition 
of the rational DAHA can
be extended to finite groups generated by complex reflections
(Dunkl, Opdam, Malle).
There is also a generalization due to Etingof- Ginzburg 
from \cite{EG} (the symplectic reflection algebras).

\rem
Following \cite{CO}, there is a one-step limiting procedure from 
$\HH $ to $\HH''.$
We set 
$$
Y_b=\exp(-\sqrt{\mathfrak{u}} D_b),\ X_b=e^{\sqrt{\mathfrak{u}} x_b},
$$ 
assuming that 
$q=e^{\mathfrak{u}}$ and tend $\mathfrak{u}\to 0.$ 
We come directly to the 
relations of the rational DAHA and the formulas for $D_b.$ 
The advantage of this direct 
construction is that the automorphisms $\tau_{\pm}$ 
obviously survive in the limit.
Indeed, $\tau_+$ in $\HH\, $ 
can be interpreted as the formal conjugation by
the $q$-Gaussian $q^{x^2/2},$ where $x^2=\sum_i x_{\om_i}
x_{\al_i^\vee}.$
In the limit, it becomes the conjugation by $e^{x^2/2},$ preserving
$w\in W,\, x_b,\, $ and taking $D_b$ to $D_b-x_b.$ Respectively, $\tau_-$
preserves $w$ and $D_b,$ and sends $x_b\mapsto x_b-D_b.$ These automorphisms
do not exist in the $\HH'.$ 
\qed
   
The {\em abstract Lusztig map} 
from $\HH'$ to $\HH''$ is as follows.
Let $w\mapsto w.$ We expand $X_\al$ in terms
of $x_\al$ in the formulas for the trigonometric 
Dunkl operators $\d_b:$ 
\begin{align}
&\d_b=\frac{1}{\mathfrak{w}} 
D_b -(\rho_k,b)+ \sum_{\al\in R_+}k_\al(b,\al)
\sum_m^\infty \frac{B_m}{m!}(-\mathfrak{w} x_\al)^m\, (1-s_\al)
\label{lusratx}
\end{align}
for the Bernoulli numbers $B_m.$ 
Then we can use them as abstract expressions for $y_b$
in terms of the generators of $\HH''.$ 

One obtains an isomorphism 
$\kapp'':\Q[[\mathfrak{w}]]\otimes\HH'\to \Q[[\mathfrak{w}]]\otimes\HH'',$
which maps $\HH'$ to the
extension of $ \HH''$  by the formal series in terms of 
$\mathfrak{w} x_b.$ An arbitrary
representation $V''$ of $\HH''\ $ which is a union of finite dimensional
$\Q_k[x]$-modules becomes an $\HH'-$ module provided that
\begin{align}
&\mathfrak{w}\ze_\al\not\in 2\pi\imath(\Z\setminus \{0\})\for
x_b(v)=(\ze,b)v, 0\neq v\in V''.
\label{intercondd}
\end{align} 

Similar to (\ref{intercond}),
this constraint simply restricts choosing $\mathfrak{w}\neq 0.$ 
The formulas for $y_b$ become locally finite
in any representations of $\HH''$ where $x_b$ act locally nilpotent,
for instance, in finite dimensional $\H''$-modules.
In this case, there are no restrictions
for $\mathfrak{w}.$ 

Finally, the composition 
$$
\kapp\equal\kapp''\circ\kapp': 
\HH[[\mathfrak{v},\mathfrak{w}]]\to 
\HH''[[\mathfrak{v},\mathfrak{w}]]
$$
is an isomorphism. Without the completion,
it makes an arbitrary finite dimensional
$\HH''$-module $V''$ a
module over $\HH$ as $q=e^{\mathfrak{v}}, t_\al=q^{k_\al}$
for sufficiently general (complex) nonzero numbers 
$\mathfrak{v},\mathfrak{w}.$ 
Note that isomorphism was discussed in \cite{BEG} (Proposition 7.1).

The finite dimensional representations are the most
natural here because, on one side,
$\kapp''$ lifts the modules which are unions of
finite dimensional $x$-modules to those for $X$, on
the other side, $\kapp'$ maps 
the $\HH'\,$-modules which are
unions of finite dimensional $y$-modules to those for $Y.$ 
So one must impose these conditions for both $x$ and $y.$
Using $\kapp$ for infinite dimensional representations
is an interesting problem. It makes the theory
analytic. For instance, the triple composition
$\kapp''\circ \g \circ \kapp'$ for the inverse Opdam transform $\g$
(see \cite{O2} and formula (6.1) from \cite{C8}) embeds $\HH\, $ in
$\HH''$ and 
identifies the $\HH''$-module $\C_c^\infty(\R^n)$ with
the $\HH$-module of PW-functions under the condition $\Re k>-1/h.$
See \cite{O2,C8} for more detail. 

\medskip
{\bf Gordon's theorem.} 
Let $k_{\sht}=-(1+1/h)=k_{\lng}, $ 
$h$ be the Coxeter number $1+(\rho,\vth).$ 
The polynomial representation $\Q[x]$ of $ \HH''$ 
is generated by $1$ and is $\{D,W\}$-spherical ($1$ is the
only $W$-invariant polynomial killed by all $D_b$). Therefore it
has a unique nonzero irreducible quotient-module.
It is of dimension $(1+h)^n,$ which was checked in \cite{BEG}, 
\cite{Go}, and also
follows from \cite{C12} via $\kapp.$

Actually, a natural setting for Gordon's theorem is with 
the parameters $k_\nu=-(e_\nu+1/h)$
for arbitrary integers $e_{\sht}, e_{\lng}\ge 0.$
The $W$-invariants and $W$-antiinvariants
of the corresponding perfect modules are connected by
the shift operators.
Cf. Conjecture 7.3 from \cite{BEG}. 
Such $k$ will not be considered in the paper.

The application of this representation to
the coinvariants of the ring of commutative 
polynomials $\Q[x,y]$ with the diagonal action of $W$ is as follows.

The polynomial representation $\Q[x]$ is naturally a quotient of 
the {\em linear space} $\Q[x,y]$
considered as an induced $\HH''$-module from the one dimensional
$W$-module $w\mapsto 1.$ So is $V''.$ The subalgebra $(\HH'')^W$ of
the $W$-invariant elements from $\HH''$ preserves $\Q\de\subset V''$
for $\de\equal\prod_{\al>0} x_\al.$

Let $I_o\subset (\H'')^W$ be the ideal of the elements vanishing
at the image of $\de$ in $V''.$ Gordon proves that 
{\em $V''$ coincides 
with the quotient $\tV''$ of $\HH''(\de)$ by the $\HH''$-submodule 
$\HH'' I_o(\de).$} It is sufficient to check that $\tV''$ is
irreducible.

The graded space gr$(V'')$ 
of $V''$ with
respect to the total $x,y$-degree of the polynomials 
is isomorphic as a linear space to the quotient of 
$\Q[x,y]\de$ by the
graded image of $\HH''I_o(\de).$  The latter
contains $\Q[x,y]^W_o\de$ for the ideal 
$\Q[x,y]^W_o\subset \Q[x,y]^W$ 
of the $W$-invariant polynomials without the constant term. 
Therefore $V''$  becomes a certain quotient of 
$\Q[x,y]/(\Q[x,y]\Q[x,y]_o^W).$ 
See \cite{Go},\cite{Hai} about the connection with
the Haiman theorem in the $A$-case and related questions for
other root systems.

The irreducibility of the $\tV''$ above is the key fact. The proof
from \cite{Go} requires considering the KZ-type local systems.
We demonstrate that the irreducibility 
can be readily proved in the $q,t$- case using the passage to
the roots of unity and therefore gives an entirely algebraic
and simple proof of Gordon's theorem via the $\kapp$-isomorphism.  

\medskip
We note that the $q,t$-generalization $V$ of $V''$ is in many ways
simpler than $V''.$ For instance, dim($V$) can be readily calculated.
However the filtration of $V''$ with respect to the degree of 
polynomials is a special feature of the rational limit as well as
the character formula from \cite{BEG, Go}, although the 
corresponding resolution has a $q,t$- counterpart.  
We will consider next the special situation when $1+h=p$ 
is prime and the definition field is $\F_p$ (see below).
In this case, the required filtration can be 
defined in the $q,t$-setting as well as for the rational
degeneration.

\medskip
{\em Modular reduction}. Concerning the same problems over $\Z$,
we may define the polynomial representation and $V''$ over $p$-adic
numbers $\Z_p$ and take its fiber over $\F_p=\Z/(p).$
There is an interesting example when $1+h$ is a prime number $p.$ 
Then $k=0$, and $\F_p\otimes\HH''$ becomes the algebra
of differential operators in $x$ with the coefficients in $\F_p.$
Its unique irreducible representation over $F_p$ with the 
nilpotent action of $x,y$ is the space 
$\F_p[x]/(x^p)$ of dimension
$p^n.$ 

This proves the theorem from \cite{Go}
for such $h.$ Indeed,
the module $\F_p\otimes \tV''$ is a semisimple $W$-module which
contains a unique one dimensional submodule with the character 
$w\mapsto \sgn(w).$ Therefore it has to coincide with $\F_p[x]/(x^p).$
Since $\F_p\otimes \tV''$ is irreducible, so is $\tV''.$

An immediate application of this construction
is that the space of diagonal coinvariants modulo $p=1+h$
has a "natural" quotient of dimension $p^n$
isomorphic to $\F_p\otimes V''.$

Actually our general prove has something in common with
this argument. However it goes via
the roots of unity instead of the modular reduction and holds for
arbitrary $h.$

\vskip 0.2cm
\section {General case} 
\setcounter{equation}{0}

The $t$-counterpart of the
element $\de$ is  
$$
\De=\prod_{\al\in R_+}
(t_\al^{1/2}X_\al^{1/2}-t_\al^{-1/2}X_\al^{-1/2}).
$$
It plays the key role in the
definition of the $t$-shift operator (see \cite{C2}). One has:
$$
T_i(\De)=-t_i^{-1/2}\De, \ 1\le i\le n.
$$
We extend $\varpi_-(T_i)=-t_i^{-1/2}$ to a one dimensional representation 
of the nonaffine Hecke algebra $\H$ generated by
$T_i, 1\le i\le n.$

Coming to a $q,t$-generalization of Gordon's 
theorem, let $q$ be generic,
$k_{\lng}=-(1+1/h)=k_{\sht}$ for the 
Coxeter number $h$. We denote the field of rationals of 
$\Q_{q,t}$ by $\tilde{\Q}_{q,t}.$

\begin{theorem}\label{GORDQT}
i) The polynomial representation $\Q_{q,t}[X]$ of $\HH$
has a unique nonzero quotient $V$ which is torsion free
and irreducible over $\tilde{\Q}_{q,t}.$ It is of dimension
$(1+h)^n.$ The action of $X$ and $Y$ is semisimple with
simple spectra. The module $V\otimes \tilde{\Q}_{q,t}$ 
considered as an $\H\,$-module contains a unique submodule 
isomorphic to $\varpi_-.$ 

ii) The module $V$ coincides with the quotient $\tV$ of 
$\Q_{q,t}[X]$ by the 
$\HH$-submodule $\HH \i_o(\tilde{\Q}_{q,t}\De)$
intersected with $\Q_{q,t}[X],$ where $\i_o$ is the kernel
of the algebra homomorphism 
$\HH_{\inv}\ni H\mapsto H(\De)\in V$ for the 
subalgebra $\HH_{\inv}$ of the 
elements of $\HH$ commuting with $T_1,\ldots,T_n.$  
\end{theorem}
 
{\em Proof.} Concerning the existence of $V$ and its isomorphism 
with the space Funct$[P/(1+h)P],$ see Theorem 8.5 and formula
(8.32) from \cite{C12}. The description there is 
for general $k=-r/h$ as $(r,h)=1.$ 
Actually we need here only the self-duality of $V,$
that is, the action of the involution $\vep$ of $\HH$ from (\ref{vepx}) in $V.$ It readily follows from the realization of $V$ 
as the quotient of $Q_{q,t}[X]$ by the radical of the
invariant bilinear form from \cite{C12}, Lemma 8.3.
In fact, the self-duality will be needed only
upon the specialization $\bullet$ below. 

We mention that
the automorphisms $\tau_{\pm}$ can be defined in $V$ as well,
but we do not use it in the paper.

The construction of $V$ in \cite{C12} holds over $\Q_{q,t}.$
Actually it suffices to have the definition of $V$ and
to prove the theorem over 
$\tilde{\Q}_{q,t}.$ Then one can use
the standard facts about the modules over PID.
Note that the module $\Q_{q,t}[X]/\HH \i_o(\De),$ generally
speaking, has torsion.
\smallskip
 
The uniqueness of $\varpi_-$ in $V\otimes \tilde{\Q}_{q,t}$ 
results from the following fact:   
\centerline{\em $W(\rho)$ is a unique 
simple $W$-orbit in $P/(1+h)P,$} 
which will be checked below.

\smallskip
\rem
There is another proof based on the shift operator, which
identifies the $\varpi_-$-component of $V=V^k$ with
the $\varpi_+$-component of $V^{k+1}$ defined for $k+1=-1/h,$ where 
$\varpi_+: T_i\mapsto t_i^{1/2}.$ The $V^{k+1}$ is one 
dimensional and coincides with its $\varpi_+$-component. The shift
operator here is the division by $\De.$ This description is
convenient to calculate the ideal $\i_o.$ Its intersections  
with $\Q_{q,t}[X]_{\inv}$ and $\Q_{q,t}[Y]_{\inv}$ are not difficult
to describe.

The self-dualty of $V$ combined with the uniqueness of $\varpi_-$
in $V$ give formally that $\tilde{V}$ is self-dual too.
Indeed, the  character $\varpi_-$ is $\vep$-invariant. So
are $\Q_{q,t}\De\subset V,$ $\i_o,$ and the kernel of the map 
$\Q_{q,t}[X]\to$ $\tV.$ The explicit discriminant
formula for $\De$ is not helpful in checking that $\vep(\De)$
is proportional to $\De$ in $V.$

The self-duality of $\tV$ can be also
deduced from the $\vep$-invariance of $\i_o,$ which can be
seen directly. \qed 
  
For $N\equal 1+h,$ we take  
$q=\exp(2\pi\imath/N)$ making 
$k=0, t^{\pm 1/2}=1.$ 
Using that $\nu_\al$ and the index of $P/Q$ are relatively prime
to $N,$ we will pick $q^{1/m}$ in the roots of unity of the 
same order
$N.$ The ${}^\bullet$ will be used to denote this specialization,
The algebra $\HH^\bullet$ is nothing else but the semidirect product
of the Weyl algebra generated by pairwise commutative
$X_a, Y_b$ for $a,b\in P$
and  $\H^\bullet=\Q_q W.$ The relations are 
$$
w\,X_a\,w^{-1}=X_{w(a)},\, w\,Y_b\,w^{-1}=Y_{w(b)},\  
X_a Y_b X_a^{-1}Y_b^{-1}=q^{-(a,b)}.
$$

We define $V^\bullet$ as a unique nonzero irreducible quotient
of $\Q_q[X].$ It is self-dual. It results from 
\cite{C12}. However it is straightforward to check 
it directly in the $\bullet$-case as well as
the semisimplicity of $X$ and $Y.$ 

Since all eigenvalues
of $Y_b$ in $V^\bullet$ (and in the whole $\Q_q[X]$) 
are $N$-th roots of unity,
the same holds for $X_b$ in $V^\bullet$ thanks
to the self-duality. Thus 
$X_b^N=1=Y_b^N$ in $V^\bullet$ for all $b\in P.$

Theorem 8.5 from \cite{C12} guarantees that $V$ remains 
irreducible under such reduction, so $V^\bullet$ is
the specialization of $V.$
It can be also seen from
the dimension formula in the following lemma.

\begin{lemma} i) The algebra $\HH^\bullet_N\equal
\HH^\bullet/(X^N=1=Y^N)$ 
has a unique irreducible
nonzero representation $V^\bullet$ up to isomorphisms.
Its dimension is $N^n.$ 

ii) It is isomorphic to $Q_q[P/NP]$ as a $W$-module. The representation 
$\varpi_-^\bullet: w\mapsto \sgn(w)$ has multiplicity one in $V^\bullet.$

iii) The quotient $\tV^\bullet$ of 
$\Q_{q}[X]$ by $\HH^\bullet \i_o^\bullet(\De^\bullet)$ is an
$\HH^\bullet_N$-module and coincides
with $V^\bullet.$
\label{WEYLQ} 
\end{lemma}

{\em Proof.} 
In the first place, $\HH^\bullet_N$ is
a group algebra of a finite group, therefore semisimple.
Let us use the well known fact
that the  Weyl algebra generated by $X_a,Y_b$
modulo the (central) relations $X^N=1=Y^N$
has a unique irreducible representation up to isomorphisms and 
$W$ acts in this representation. This representation
equals $\Q_q[X]/(X^N=1)$ and  
is nothing else but $V^\bullet.$ 

The multiplicity one statement from ii)
follows from the uniqueness of a simple $W$-orbit in $P/NP.$
Let us check it. 

We can assume that the orbit is $W(b)$
for $b\in P$ such that $0< (b,\al^\vee)< 1+h$ for all $\al\in R_+.$
Therefore $b$ can be $\rho=\sum_i^n \om_i$ or 
$\rho+\om_r $ for $r\in O',$
i.e., for a minuscule weight $\om_r.$ Indeed, the coefficient of
$\al_i^\vee$
in the decomposition of $\vth$ in terms of simple coroots is one
only for $r\in O'.$ For $b=\rho+\om_r,$ let $w=u_r^{-1}$
for $u_r$ from $\om_r=\pi_r u_r.$ Then
\begin{align}
&(w(\rho)-b,\al_r)= -(\rho,\vth)-(\om_r,\al_r^\vee)-1=-(1+h)=-N\and
\notag\\ 
&(w(\rho)-b,\al_i)=(\rho,\al_j)-(\rho,\al_i)=0 \for i\neq r,\
 w^{-1}(\al_i)=\al_j.\notag
\end{align}
Thus $b$ and $\rho$ generate the same $W$-orbit modulo $NP.$

The module $V^\bullet$ is self-dual. So is $\tV^\bullet$ because
the $\H$-module $\Q_q \De^\bullet$ is of multiplicity one in 
$\tV^\bullet$ and therefore invariant with respect to $\vep.$ 
It gives that $\tV^\bullet$ is a finite dimensional 
$\HH^\bullet_N$-module. It  has $V^\bullet$ as a quotient, and contains a unique $W$-submodule isomorphic to $\varpi_-^\bullet.$ 

Supposing that the kernel $K$ of the map $\tV^\bullet\to 
V^\bullet$ is nonzero,
it must contain a nonzero $\HH^\bullet_N$-submodule.  
Hence $K$ contains at least one copy of $V^\bullet$
(the uniqueness), and the multiplicity of $\varpi_-^\bullet$ in 
$\tV^\bullet$ cannot be one.
\qed 

We would like to mention that claims
(ii)--(iii) are somewhat unusual in the general theory of
Weyl algebras. Given the rank,
they hold only for special choice of roots of unity.
For instance, $N$ must be $3$ in the case of $A_1.$
\medskip

Coming back to the general case,
the coincidence statement
of the theorem is actually over $\tilde{\Q}_{q,t},$
so it suffices to check it at one special point.
The lemma gives it for the $\bullet$-point. To be more exact, 
claim (iii) of the lemma gives the irreducibility of $\tV$
and therefore the coincidence $\tV=V$ at the common point and
the theorem.   
\qed

\smallskip
The applications to the diagonal coinvariants goes via the
universal Dunkl operators, which will be introduced in the
next section.  The problem is to calculate the action 
of $Y_b$ in the linear space of Laurent polynomials
$\Q_{q,t}[X,Y]$ identified with the $\HH\,$-module induced from
a one dimensional character of $\H.$ We assume here 
that the $X$-monomials
are placed before $Y$-monomials. Then the action of $X_b$ in 
the space $\Q_{q,t}[X,Y]$ is the "commutative" 
multiplication by $X_b.$ The action of $Y_b$
is by the left multiplication of the monomials in the form
$X^{\cdot} Y^{\cdot}.$ Hence it requires reordering and
leads to nontrivial formulas. 

The $\HH\,$-module 
$\Q_{q,t}[X,Y]$ is obviously self-dual.
However since it is necessary to order $X$ and $Y$ after applying
$\vep,$ the formulas for its action in $\Q_{q,t}[X,Y]$ 
are involved.
\smallskip
 
Now, the module $V$
is the quotient of the module $\Q_{q,t}^-[X,Y]$ induced for the 
character $\varpi_-$ of $\H$ by its
submodule $\HH \i_o (1).$ It is a one-parametric
deformation of
the polynomial ring $\Q_{q,t}[X,Y]$ divided by the 
ideal $\Q_{q,t}[X,Y]_o^W=$ 
$\{g(X,Y)(f(X,Y)-f(1,1))\}$ for $W$-invariant Laurent polynomials
$f$ and arbitrary $g.$ Therefore it can be identified
with a quotient
of the space of diagonal coinvariants 
$$
\Q[X,Y]/(\Q[X,Y]\Q[X,Y]_o^W)\simeq
\Q[x,y]/(\Q[x,y]\Q[x,y]_o^W),
$$
with the ring of definition extended to $\Q_{q,t}.$

Concerning Gordon's theorem,
the degenerations $V',V''$ of $V$ can be introduced 
as quotients of the polynomial representation
by the radicals of the degenerations of
the bilinear form from \cite{C12}, Lemma 8.3.
They are irreducible modules
for $\HH',\HH''$ thanks to Lusztig's isomorphisms, and
satisfy the same multiplicity one statement. The modules
$\tV',\tV''$ are defined in terms of the ideal of the
$W$-invariant elements of the double Hecke algebra
vanishing at the discriminant 
subspace of $V',V''$ (corresponding to the sign-character of $W$). 
The $\tV',\tV''$ are irreducible due to Lusztig's isomorphisms,
and therefore $V'=\tV'$ and 
$V''=\tV''.$ The latter is
the (main part of the) Gordon's theorem.
   
One can also consider the specialization
$\HH^\bullet$  assuming that $N$ is a 
prime number $p$ and make the field of constants  
$\F_p.$ Then the space $\F_p[X,Y]/(\F_p[X,Y]\F_p[X,Y]_o^W)$
has a "natural" quotient-space isomorphic to $V^\bullet$ over $\F_p.$
In this case, the formulas for the Lusztig
homomorphism $\kapp$ contain no denominators divisible by $p$
and it is well defined over $\F_p.$ Applying
$\kapp$, we establish the coincidence of this quotient with the one
obtained at the end of Section \ref{SECTDEGEN} using the 
differential operators over $\F_p.$

\vskip 0.2cm
\section{Universal DAHA}
First, we will give a   
$X\leftrightarrow Y$--symmetric presentation
of $\HH\ .$ It goes via the {\em universal double affine braid group} 
$\widehat{\mathfrak{B}}.$ This group is defined
to be generated by the pairwise commutative
$X_b$, the pairwise commutative  $Y_b,$ the elements 
$\widehat{T}_i,$ where $b\in P,$ $0\le i\le n,$  and the group 
$\widehat{\Pi}=\{\widehat{\pi}_r,r\in O\}\simeq \Pi$   
with the following defining relations:

\par (a) $ \widehat{T}_i\widehat{T}_j\widehat{T}_i... = 
\widehat{T}_j\widehat{T}_i\widehat{T}_j...,\ m_{ij}$ 
factors on each side, 
\par\ \ \ \ \,$\widehat{\pi}_r\widehat{T}_i\widehat{\pi}_r^{-1} = 
\widehat{T}_j$ if $\widehat{\pi}_r(\al_i)=\al_j$;  

\par (b) $\widehat{T}_iX_b \widehat{T}_i\ =
\ X_b X_{\al_i}^{-1},\ 
\widehat{T}_i^{-1}Y_b \widehat{T}_i^{-1}\ =\ 
Y_b Y_{\al_i}^{-1}$ 
\par\ \ \ \ \, if  $(b,\al^\vee_i)=1 \for 0 \le i\le  n$;

\par (c) $\widehat{T}_iX_b\ =\ X_b \widehat{T}_i,\ 
\widehat{T}_iY_b\ =\ Y_b \widehat{T}_i$ 
\par\ \ \ \ \, if $(b,\al^\vee_i)=0 \for 0 \le i\le  n$;

\par (d) $\widehat{\pi}_rX_b \widehat{\pi}_r^{-1} = 
X_{\widehat{\pi}_r(b)},\ 
\widehat{\pi}_rY_b \widehat{\pi}_r^{-1} = 
Y_{\widehat{\pi}_r(b)},\ r\in O'$.

\medskip
Note that no relations between $X$ and $Y$ are imposed.
We continue using the notation $X_{[b,j]}=X_b q^j,$
$Y_{[b,j]}=Y_b q^{-j}.$ The element $q^{1/m}$ is treated
as a generator which is central. Later, $q^{1/(2m)}$ will
be needed.
 
The relations (a-d) are obviously invariant with respect
to the involution  
\begin{align}\label{eptilde}
&\widehat{\vep}:\, X_b\leftrightarrow Y_b,\,
\widehat{T}_i\mapsto \widehat{T}_i^{-1} (0\le i\le n),\,
\widehat{\pi}_r\mapsto \widehat{\pi}_r (r\in O).
\end{align}
Concerning (d), recall that $\pi_r^{-1}$ is $\pi_{r^*}$
for $r^*\in O.$ The same holds for 
$u_r=\pi_r^{-1}\om_r$ and
$\widehat{\pi_r}.$

\begin{theorem}
The group $\mathfrak{B}$ generated by $X,T,\Pi,q^{1/(2m)}$ subject
to the relations (i--v) from Definition \ref{doublex}
coincides with the quotient
of $\widehat{\mathfrak{B}}$ by the relations:
\begin{align}\label{bbrela}
&\widehat{T}_0=q^{-1}\,X_\vth 
\widehat{T}_{s_\vth} Y_\vth^{-1},\ 
\widehat{\pi}_r=q^{(\om_r,\om_r)/2}\,Y_r
\widehat{T}_{u_r}^{-1}X_{r^*}^{-1}.
\end{align}
The map is 
\begin{align}\label{prmap}
\mathfrak{pr}:\, X_b&\mapsto X_b,\, Y_b\mapsto Y_b,\
\widehat{T}_i\mapsto T_i\, (i>0),\,
\widehat{T}_0\mapsto 
q^{-1}\,X_\vth T_{s_\vth} Y_\vth^{-1},\notag\\
\widehat{\pi}_r&\mapsto 
q^{(\om_r,\om_r)/2}\,Y_r T_{u_r}^{-1}X_{r^*}^{-1},\
q^{1/(2m)}\mapsto q^{1/(2m)},
\end{align}
where the elements $Y_b\in \HH\,$ are given by
(\ref{Ybx}). The images of $\widehat{T}_i,\widehat{\pi}_r$
in $\mathfrak{B}$ coincide with
$\tau_+(T_i),\tau_+(\pi_r)$   
for $\tau_+$ from (\ref{taux}). The relations
(\ref{bbrela}) are invariant with respect to the involution
$\widehat{\vep}.$ The latter becomes $\vep$ from (\ref{vepx})
in $\mathfrak{B}.$
\end{theorem}

{\it Proof.}  The key fact here is that 
\begin{align}
&\vep(\tau_+(T_i))=(\tau_+(T_i))^{-1}\, (i\ge 0),\ 
\vep(\tau_+(\pi_r))=\tau_+(\pi_r),  
\label{tautpi}
\\
&\where \tau_+(T_i)\,=\,T_i \for i>0, \ 
\tau_+(T_0)\,=\,q^{-1}\,X_\vth T_{s_\vth} Y_\vth^{-1},\,
\notag\\ 
&\tau_+(\pi_r)=q^{\frac{(\om_r,\om_r)}{2}}\,
Y_r T_{u_r}^{-1}X_{r^*}^{-1},\ r\in O.
\notag
\end{align}
See formula (2.17) for the action of the involution 
$\eta=\vep\tau_-^{-1}\tau_+\tau_-^{-1}$ from \cite{C12}.
We note that (2.17) there
directly results in the invariance of the Gaussians with
respect to the difference Fourier transform corresponding
to the involution $\vep.$ The elements 
$\tau_+(T_i),\tau_+(\pi_r)$ are exactly the images of the elements
$\widehat{T}_i,\widehat{\pi}_r$ in $\mathfrak{B}$ under 
$\mathfrak{pr}.$
Relations (\ref{tautpi}) readily follow from the explicit
formulas. 
\qed

\smallskip
There are important quotients of the group $\widehat{\BB}$ and
the algebra $\widehat{\HH}$ (see below) 
obtained by imposing the commutativity of $X$ with $Y.$
They are essential in the theory of the difference Fourier transform.
The $\{X_i\}$ are treated as the coordinates,
$\{Y_i\}$ play the role of spectral parameters. 
Note the immediate projection of these quotients to $\BB$, $\HH\,$
when we make $Y_b=X_b^{-1}.$  
\smallskip

\rem
The realization of $\mathfrak{B}$ as a quotient of 
$\widehat{\mathfrak{B}}$ can be used for a direct proof
that $\vep$ can be extended to an involution of $\mathfrak{B}.$ 
It can simplify the 
straightforward proof due to Macdonald and the
author from \cite{Ma3}, but the difference is not very significant.
Our $\widehat{\mathfrak{B}}$ has something in common with the
"triple affine Artin group" introduced recently by Ion and 
Sahi in \cite{IS} for the purpose of interpreting
the projective action of $PSL_2(\Z).$ 
Compare our relation (\ref{bbrela})
and formula (23) in \cite{IS} which establishes the
connection of their group with the double affine braid
group from \cite{C15}.
\qed

\medskip

Turning to the Hecke algebras, let us define the 
{\em universal affine double Hecke algebra} 
$\widehat{\HH}$ as
the quotient of the group algebra $\Q_{q,t}\widehat{\mathfrak{B}}$
by the quadratic relations (o) from Definition \ref{doublex}
for $\widehat{T}_i.$ Here we do not assume that 
$t_{\lng}=t_{\sht}.$
Recall that the elements $X_b$ and $Y_b$ are entirely
independent, so the counterpart of the PBW theorem is
that  an arbitrary element $\widehat{H}\in 
\widehat{\HH}$ can be uniquely
represented as 
$$
\widehat{H}=\sum_{\hw\in \hW} Q_{\hw}T_{\hw},\ \where Q_{\hw} 
\hbox{\ are\ noncommutative\ polynomials\ in\ } X,Y. 
$$

For applications to the Dunkl operators,
this definition will be needed in the following 
form. We claim that the algebra $\widehat{\HH}$ is
generated over $\Q_{q,t}$ by the affine Hecke algebra 
$$
\widehat{\h}\equal
\lan \widehat{T}_i,\widehat{\pi}_r\ran, \ i\ge 0,\ r\in O,
$$
the pairwise commutative $X_b\, (b\in P)\, ,$ the pairwise 
commutative $Y_b\, (b\in P),$ satisfying the relations
(d) above with the $\widehat{\pi}_r,$
and the Lusztig-type relations
\begin{align}
&\widehat{T}_iX_b -X_{s_i(b)}\widehat{T}_i\ =\ 
(t_i^{1/2}-t_i^{-1/2})\frac{X_{s_i(b)}-X_b}{
X_{\al_i}-1},\ 0 \le i\le  n,
\label{tixix}\\
&\widehat{T}_iY_b -Y_{s_i(b)}\widehat{T}_i\ =\ 
(t_i^{1/2}-t_i^{-1/2})\frac{Y_{s_i(b)}-Y_b}{
Y_{\al_i}^{-1}-1},\ 0 \le i\le  n.
\label{tiyix}
\end{align}
Imposing (\ref{bbrela}), we represent $\HH\,$ as a quotient
of $\widehat{\HH}.$ Note that
this definition is compatible with the restriction
to the lattices $B$ between $Q$ and $P$ taken instead of $P.$ 

We set
$$
\widehat{Y}_{b_+}\equal \widehat{\pi}_r
\widehat{T}_{i_1}\cdots\widehat{T}_{i_l} \for
b_+=\pi_r s_{i_1}\cdots s_{i_l} \hbox{\ as\ } b_+\in P_+,
\,l=l(b_+),
$$
more generally, $\widehat{Y}_{b_+-c_+}=\widehat{Y}_{b_+}
\widehat{Y}_{c_+}^{-1}$ for $b_+,c_+\in P_+.$
\medskip

{\em Universal Dunkl operators.}  
Given a representation $\widehat{V}$ of $\widehat{\h},$
the general  universal Dunkl operators are the images of
$\widehat{Y}_b \, (b\in P)\, $ 
and $\widehat{\pi}_r\, (r\in O)\, $
in the $\widehat{\HH}\,$-module $\i_{\widehat{V}}$
induced from $\widehat{V}.$ As a linear space,
it is isomorphic to the linear space of {\em noncommutative}
polynomials of $X,Y$ with the (right) coefficients in 
$\widehat{V}\,:\,$
$\i_{\widehat{V}}=\cup_{e\in \N}\,\p_e$ for  
\begin{align}\label{pesum}
&\p_e\equal\bigl\{ 
\sum_{\mathbf{b},\mathbf{c}\in \mathbf{P}} 
X_{b_1}Y_{c_1}\cdots X_{b_e}Y_{c_e}
\widehat{v}_{\mathbf{b},\mathbf{c}}\bigr\},\ 
\mathbf{b}\in \mathbf{P}=P^e \ni
\mathbf{c},\, \widehat{v}_{\mathbf{b},\mathbf{c}}\in \widehat{V}. 
\end{align}

Here the sums in (\ref{pesum}) represent different
vectors in $\i_{\widehat{V}}$ for different $v$-coefficients 
if we assume that $b_i\neq 0\neq c_j$ 
for the indices $1<i\le e,$ $1\le j<e$
as $\widehat{v}_{\mathbf{b},\mathbf{c}}\neq 0.$

The action of $X$ and $Y$ is by the left multiplication.
The subspaces $\p_e$ are $\widehat{\h}$-submodules and
also $X$-submodules for $e>0.$ 

The action of $\widehat{T}_i\, (i\ge 0)\,$ and
$\widehat{\pi}_r$ in $\p_e$
can be calculated using (\ref{tixix}),(\ref{tiyix}), and
the relations (d) above.

In the first interesting case $e=1,$ it is as follows: 
\begin{align}\label{dunxy}
&\widehat{T}_i(X_bY_c\widehat{v})\ =\ 
 X_{s_i(b)}Y_{s_i(c)}\widehat{T}_i(\widehat{v})+ 
\\ 
&(t_i^{1/2}-t_i^{-1/2})
\bigl(     \frac{X_{s_i(b)}-X_b}{X_{\al_i}-1}Y_c +
 X_{s_i(b)}\frac{Y_{s_i(c)}-Y_c}{Y_{\al_i}^{-1}-1}
\bigr)\, \widehat{v},\ 0 \le i\le  n,
\notag\\
&\widehat{\pi}_r(X_bY_c\widehat v)\ = 
\ X_{\pi_r(b)}Y_{\pi_r(c)}\widehat{\pi}_r(\widehat v)\ 
\for b,c\in P, \ \widehat{v}\in \widehat{V}.
\end{align}

When the initial representation $\widehat{V}$
is the character
$$
\widehat{\varpi}_+: \widehat{T}_i\mapsto t_i^{1/2},\ i\ge 0,\
\widehat{\pi}_r\mapsto 1,
$$
we come to a double variant of
formulas (\ref{Demazx}) and the corresponding 
Dunkl operators. Namely:
\begin{align}\label{dunxyo}
&\widehat{T}_i\ =\ t_i^{1/2}s_i^xs_i^y+  
(t_i^{1/2}-t_i^{-1/2})
\bigl(\frac {s_i^x-1}{X_{\al_i}-1}+ 
 s_i^x \frac{s_i^y-1}{Y_{\al_i}^{-1}-1}\bigr). 
\end{align}
Here $s_i^x, s_i^y$ 
act respectively on $X$ and $Y,$
the differences are applied before the division in 
the divided differences. Similarly,
$\widehat{\pi}_r=\pi_r^x\pi_r^y.$

\medskip
{\em Double polynomials.}
The above consideration leads to the formulas
for the action of the hat-operators in the
$\HH\,$-module $\Q_{q,t}[X,Y]$ induced from the 
character $\varpi_+$ of the nonaffine Hecke 
subalgebra $\H.$

These operators are the images
of $\{\widehat{T}_i,\widehat{\pi}_r\}$
under the projection $\mathfrak{pr}.$ They coincide
with $\tau_+(T_i),\tau_+(\pi_r).$ See (\ref{tautpi}).
We will use the same hat-notation for them, although
now they are the elements of $\HH\, .$
 
The formulas are:
\begin{align}\label{dundub}
&\widehat{T}_i(X_bY_c)= 
 X_{s_i(b)}Y_{s_i(c)}\widehat{T}_i(1)+\notag\\
&(t_i^{1/2}-t_i^{-1/2})
\bigl(     \frac{X_{s_i(b)}-X_b}{X_{\al_i}-1}Y_c +
 X_{s_i(b)}\frac{Y_{s_i(c)}-Y_c}{Y_{\al_i}^{-1}-1}
\bigr),\notag\\
&\widehat{T}_i(1)=t_i^{1/2}\, (i>0),\ 
\widehat{T}_0(1)=q^{-1}\,
X_{\vth}T_{s_\vth}Y_\vth^{-1}(1)\notag\\
&=q^{-1}\,X_{\vth}
(Y_\vth T_{s_\vth}^{-1}(1)-(t_0^{1/2}-t_0^{-1/2})),
\notag\\
& T_{s_\vth}^{-1}(1)=t_{\sht}^{1-(\vth,\,\sum_{\sht} \al^\vee)}
\, t_{\lng}^{-(\vth,\,\sum_{\lng} \al^\vee)},\ \al\in R_+.
\end{align}

The $X$-monomials act by the left multiplication, the
operators $\widehat{\pi}_r$ via the relations (d)
and the formula for $\widehat{\pi}_r(1)$ similar to
that from (\ref{dundub}). Knowing the action of
$X_b,\widehat{T}_i,$ and $\widehat{\pi}_r$ is sufficient
for determining the structure of the $\HH\, $-module.

Recall that the involution $\vep$ acts naturally in $\Q_{q,t}[X,Y]$
sending $\widehat{T}_i\mapsto$ $\widehat{T}_i^{-1}$ (all $i$) and
$\widehat{\pi}_r\mapsto$ $\widehat{\pi}_r.$ 

We note that
the hat-operators can be used, for instance,
to introduce the "double radial parts". As usual, the
simplest ones corresponding to the miniscule symmetric monomial
functions, can be calculated explicitely.

\medskip
To conclude,
we note that there is an interesting group of automorphisms
of $\widehat{\BB},\widehat{\HH}\,$ generated by 
the tau-automorphisms
$$
\tau^x_{\pm},\ \ \tau^y_{\pm}\equal
\widehat{\vep}\tau^x_{\mp}\widehat{\vep}.
$$  
They act respectively in $\lan \widehat{T},\widehat{\pi},X\ran$  fixing $Y,$ 
and in 
$\lan \widehat{T},\widehat{\pi},Y\ran$ fixing $X.$  
This group is an extension of the
$PSL_3(\Z).$  Hopefully it is "naturally" connected with the 
action of $PSL_3(\Z)$ on solutions of some KZB-type equations 
found by Felder and Varchenko. 

\medskip

\bibliographystyle{unsrt}

\end{document}